\documentclass[a4paper,12pt]{article} 

\usepackage[T2A]{fontenc}
\usepackage[cp1251]{inputenc}
\usepackage[english]{babel}
\usepackage{amssymb,amsfonts,amsmath,mathtext,mathtools,cite,enumerate,float} 
\usepackage{graphicx} 
\usepackage{epstopdf}
\makeatletter

\newcommand{\at}[2][]{#1|_{#2}}

\renewcommand{\i}{i}
\newcommand{\ptl}{\partial}
\newcommand{\vph}{\varphi}
\newcommand{\eps}{\varepsilon}
\makeatother

\begin{document}

\title{Diffraction by an elongated body of revolution. A boundary integral 
equation based on the parabolic equation}
\author{A. V. Shanin, A. I. Korolkov}
\maketitle
% REQUIRED
\begin{abstract}
A problem of diffraction by an elongated body of revolution is studied. The incident wave falls along the axis. The wavelength is small comparatively to the dimensions of the body. 
The parabolic equation of the diffraction theory is used to describe the diffraction 
process. A boundary integral equation is derived. The integral equation 
is solved analytically and by iterations for diffraction by a cone. 
\end{abstract}

\section{Introduction}
We consider a problem of diffraction by an obstacle having geometrical sizes 
large comparatively to the wavelength.
The body is assumed to be elongated, i.~e.\ the axial (longitudinal) size
should be bigger than the radial (transversal) size. The incident wave falls along the axial direction or closely to it. These features enable one to make some 
simplifications to the otherwise general diffraction problem and to solve 
this problem \cite{Babich1979,Kirpichnikova2013, Popov2009,Popov2014e}.  
One can see that the considered problem should have three small parameters 
(related to the wavelength, transversal / longitudinal size, and the angle of incidence). These small parameters and the relations between them are discussed in the paper. 

Diffraction by elongated bodies attracts considerable attention of researchers. 
Several different approaches can be marked. First, there is a traditional asymptotic approach when the field is represented as rays. 
Within this approach the creeping  waves are generated in the Fock's zones \cite{Fock1965} and then propagate along convex surfaces or detach from the surfaces \cite{Engineer1998}.
Second, there is an approach developed in \cite{Andronov2011e,Andronov2012ae,Andronov2014e, Andronov2015e}. The field is expanded as Fourier series with respect to the rotational angle. A parabolic equation is solved or each Fourier component  in the radial and axial variables. 
Third, for conical obstacles (with not necessarily circular cross-section) one can 
apply methods which are specifically conical. They are based either on Smyshlyaev's 
formula \cite{Babich2000,Shanin2012c} or on Kontorovich-Lebedev integral representation \cite{Lyalinov2007, Antipov2002}. 

All these methods are mathematically complicated, and none of them is physically clear. Here we are trying to develop a method that could admit a simple physical interpretation. We take the parabolic equation of diffraction theory (PETD) \cite{Levy2000} as the governing equation. On its base we derive a boundary integral equation similar to the Hong's one \cite{Hong1967}. This equation is of Volterra type, so it can be solved by iterations. 

The integral equation is derived in two forms. First, this is an equation, whose kernel depends on four degrees of freedom: two coordinates of the source
and two coordinates of the receiver (both belong to the surface of the obstacle). Such an equation can be written for any elongated body, not 
necessarily for a body of revolution. Then, using  the rotational symmetry, a 1D integral equations for each angular mode is derived. If the axial incidence is considered, the consideration can be  restricted to the equation for the zeroth mode. 

The integral equation can be solved analytically or numerically. Also it is important that the 
equations belong to the Volterra class, so they can be solved by iterations. 
We demonstrate the capabilities of the method on the example of diffraction by a cone.   

\section{Problem formulation}

\subsection{Problem formulation for the Helmholtz equation}

Consider a 3D space, in which the cylindrical coordinates  
$(x,r, \varphi)$ are set. The $x$ direction is assumed to be {\em longitudinal}.
We assume that all wave processes are paraxial, i.~e.\ the most
part of angular spectrum is concentrated near the positive $x$ direction.  

The Helmholtz equation 
\begin{equation}
\label{helm}
\left(\frac{\ptl^2}{\ptl x^2}+ \Delta_\perp + k^2\right)\tilde u(x,r,\varphi)=0,
\end{equation}
\[
\Delta_\perp \equiv
\frac{1}{r}\frac{\ptl}{\ptl r}r\frac{\ptl}{\ptl r}+ \frac{1}{r^2}\frac{\ptl^2}{\ptl \varphi^2}
\]
is fulfilled in the space outside the scattering body. 
The time dependence is 
assumed in the form $e^{-i \omega t}$ and is omitted henceforth. 
The body occupies the domain 
$r < f(x)$ having the rotational symmetry. 
The body surface $r = f(x)$ will be denoted by $\Gamma$.
The axis of the body can be compact ($X_1 \le x \le X_2$) or half-infinite
($X_1 < x$). An important example of the body is 
a cone: 
\begin{equation}
r = \alpha x, \qquad x > 0.
\label{coneshape}
\end{equation}

The incident wave has form
\begin{equation}
\label{tuin}
\tilde u^{\rm in} = \exp\{i k(x \cos \theta + r \cos \vph \sin \theta \},
\end{equation}
where $\theta$ is the angle of incidence. We assume that $\theta$ is small. Moreover, in 
all examples we assume that $\theta = 0$, i.~e.\ the incidence is axial. 

Total field $\tilde u$ is represented as follows:
\begin{equation}
\tilde u = \tilde u^{\rm in} + \tilde u^{\rm sc} .
\end{equation}
Here $\tilde u^{\rm sc}$ is the scattered field. The scattered field should satisfy the  radiation condition. We formulate the radiation condition in the form of the limiting absorption 
principle. Namely, we assume that $k$ has a small positive part. Then the scattered part should decay as $|x|\to \infty$ and $r\to \infty$. 

Neumann boundary condition $\ptl \tilde u / \ptl n =0$ is fulfilled on the surface of the body. Here $n$ is the outward normal to the boundary. 

If the body shape has conical points, 
the total field should satisfy Meixner's condition at these points. 
The restriction is that the ``energy'' combination
$|\nabla \tilde u|^2 + \tilde u^2$ should be locally integrable near the tip. 

%%%%%%%%%%%%%%%%%%%%%%%%%%%%%%%%%%%%%%%%%%%%%%%%%%%%%%%%%%%%%%%%%%%%%%%%%%%%%%%%%%%%%%%%%%%%

\subsection{Parabolic formulation}

The body of revolution is assumed to be {\em elongated}. 
We are planning to use  the
parabolic equation, so most of the waves should be scattered on small angles.  
It is not, however, simple to write down a general definition of an elongated body. 

The main example here is a cone (\ref{coneshape}), for which  
the angle $\alpha$ and the incidence angle $\theta$ should be small comparatively to~1. 

A formal asymptotic description of an elongated body can be found in 
\cite{Popov2009}. Two curvatures of the shape are introduced at each point of the body, the 
transversal one $K_t = (f(x))^{-1}$, and the longitudinal one, $K_l = \dot f(x)$ 
(here and below, the dot denotes the derivative with respect to the argument).
The wavenumber $k$ tends to infinity. The authors of \cite{Popov2009} build an asymptotic expansion 
of the field with respect to two large parameters: 
\[
M_0 = (k/2 K_l)^{1/3}
\mbox{ and }
\Lambda_0 = k K_t.
\]    
The expansion built in \cite{Popov2009} is valid when the large parameters are linked via the relation
$\Lambda_0 = M_0^{2- \nu}$, where $0<\nu<2$. This means that not only $k\to \infty$, but also 
that the shape of the body evolves in the asymptotic process. 

In \cite{Andronov2012ae} the case $\Lambda_0 = M_0^2$ is studied. In this case one cannot describe the field in the ray form 
and should consider the field not locally with respect to $\vph$.

Here we are not trying to build a full asymptotic expansion of the field. Instead, we just control that the scattering process is paraxial, i.~e.\ that the angle of scattering is small. 
Three conditions should be valid: 

--- The angle of incidence should be small: $\theta \ll 1$.

--- The slope of the surface  should be small: $\dot f \ll 1$. 

--- The angle of diffraction by a convex surface should be small: $(\ddot f / k)^{1/3} \ll 1$.

The last expression can be commented. The longitudinal size of the Fock's zone 
is $\Delta x \sim  (\ddot f)^{-2/3} k^{-1/3}$ \cite{Fock1965}. 
The diffraction angle can be estimated as the change of the surface slope   
on this size, i.~e.\ as $\Delta x \ddot f = (\ddot f / k)^{1/3}$.

Under the listed conditions, 
represent the field as
\begin{equation}
\label{upar}
\tilde u(x,r,\varphi)= \exp\{i k x\}u(x,r,\varphi),
\end{equation}
$u$ is a slow function of $x$, as compared to the exponential factor.
Traditionally, $u$ is called the attenuation function.
 Substituting (\ref{upar}) in (\ref{helm}) and neglecting the term with the second $x$-derivative (see \cite{Vlasov1995}), get an approximate equation for $u$:
\begin{equation}
\label{pareq}
\left(\frac{\ptl}{\ptl x}+ \frac{1}{2ik} \Delta_\perp \right)u = 0,
\end{equation}
This is the parabolic equation of the diffraction theory (PETD).

The incident wave (\ref{tuin}) in the parabolic approximation looks as follows:
\begin{equation}
\label{uin1}
u^{\rm in}(x,r,\vph) = \exp\left\{i k (\theta r \cos\vph - x \theta^2 /2)\right\}.
\end{equation}
This expression takes into account that $\cos \theta \approx 1- \theta^2 / 2$, 
$\sin \theta \approx \theta$. Moreover, it is easy to check that (\ref{uin1}) obeys
(\ref{pareq}).

In the case of axial incidence $\theta = 0$, and
\begin{equation}
\label{uin}
u^{\rm in} =1.
\end{equation} 

Neumann boundary condition takes the following form:
\begin{equation}
\label{boundcond}
N [u](x, f(x),\vph) =0,
\qquad 
N \equiv \frac{\ptl u}{\ptl r} - ik \dot f. 
\end{equation}
This boundary condition can be explained as follows. 
Consider the normal derivative of the field $\tilde u$:
\[
\frac{\ptl \tilde u}{\ptl n} = 
\frac{e^{ikx}}{\sqrt{1+(\dot f)^2}} \left( 
\frac{\ptl u}{\ptl r} - \dot f \frac{\ptl u }{\ptl x} - i k \dot f u 
\right)  .
\]
Approximate $\sqrt{1+(\dot f)^2}$ as 1 and omit the second term 
(it is smaller than the third one). As the result get (\ref{boundcond}).  
A detailed procedure of derivation of this boundary condition 
and a discussion of its properties can be found in \cite{Shanin2016en}. 

Condition (\ref{boundcond}) can be rewritten in the following way:
\begin{equation}
\label{boundsc}
N[u^{\rm sc}](x,f(x),\varphi) = -N[u^{\rm in}](x,f(x),\varphi).
\end{equation}

The scattering problem for the parabolic equation has the nature quite 
different from the Helmholtz one. 
This equation describes only waves traveling in the positive direction. 
Thus, if there is no obstacle at $x < X_1$, there should be no scattered waves as well. 
Therefore, we should formulate the initial condition 
\begin{equation}
u^{\rm sc} (x, r, \vph) = 0 
\quad \mbox{ for } x < X_1.
\label{initcond}
\end{equation}

A radiation condition in the limiting absorption form would be as follows. For $k$ having 
a small positive imaginary part the scattered field should decay as $r \to \infty$. 

Finally, the problem should include a ``Meixner's'' condition for the conical points. 
The condition has the following form. The combination 
\[
\left| \nabla_\perp u \right|^2 + |u|^2 
\]
should be integrable in any small domain of the space outside the scatterer. 
Here $\nabla_\perp$ is the gradient in the plane normal to the $x$-axis.
This is an analog of the energy-like combination for the Helmholtz equation.

The ``Meixner's'' condition provides uniqueness of the solution, however the solution of the 
parabolic equation near the conical point is not close to the solution of the Helmholtz equation. The zone in which the solutions are different in the $x$-direction is about several 
wavelengths (i.~e.~$\sim k^{-1}$).

The parabolic equation studied in this paper is {\em global\/}, i.~e.\ it is valid everywhere outside the obstacle. This contrasts with the approach used by authors of \cite{Popov2009,Andronov2012ae}, 
who use the parabolic equation only in a narrow boundary layer near the surface of the obstacle. 
However, this difference is not important since anyway only the field on the surface is needed to reconstruct the directivity diagram.

The main benefit of using the parabolic equation (comparatively to the Helmholtz equation) is that the spatial scale of the wave process becomes of order of the Fock's zone, i.~e.\ it becomes much bigger than the wavelength. If, say, a numerical method is used to solve the problem, 
one should take about 10 nodes per the size of the Fock's zone rather than 10 points per a wavelength, so the numerical procedure demands less resources.

%%%%%%%%%%%%%%%%%%%%%%%%%%%%%%%%%%%%%%%%%%%%%%%%%%%%%%%%%%%%%%%%%%%%%%%%%%%%%%%%%%%%%%%%%

\section{Derivation of boundary integral equations}

\subsection{Green's function and Green's theorem for the parabolic equation}

Write down the Green's function of the parabolic equation. 
Denote points of the space by straight letters ${\rm r} = (x, r, \vph)$. 
Let a point source be located at ${\rm r }_s = (x_s , r_s , \vph_s)$.
Introduce the Green's function as a solution of an inhomogeneous 
parabolic equation
\begin{equation}
\left( \frac{\ptl}{\ptl x} + \frac{1}{2i k }
\Delta_\perp \right)
g({\rm r}, {\rm r}_s) = \delta({\rm r} - {\rm r}_s),
\label{Greens}
\end{equation}
where the operator in the left acts on the components of ${\rm r}$, $\delta$
is the Dirac's delta-function. The solution should obey the 
initial condition, i.~e.\ it should be equal to zero if $x < x_s$.  

One can check directly that 
the Green's function has the following form for $x > x_s$:
\begin{equation}
\label{greensfunction}
g({\rm r}, {\rm r}_s) = \frac{k}{2 \pi i (x - x_s)}
\exp \left\{ \frac{i k}{2} \frac{(\Delta r)^2}{x- x_s}\right\},
\end{equation}
$\Delta r$ is the distance between the projections of ${\rm r}$ and ${\rm r}_s$
onto the transversal plane:
\[
(\Delta r)^2 = r^2 + r_s^2 - 2 r r_s \cos(\vph - \vph_s).
\]

An important property of the Green's function is as follows. Let $v$ be a solution of the 
homogeneous parabolic equation (\ref{pareq}) everywhere in the 3D space, for example 
$v$ can be equal to $u^{\rm in}$. Then  
for 
\[
{\rm r} = (x, r, \vph), \qquad {\rm r}_* = (x_*, r_*, \vph_*), \qquad  x_* > x,
\] 
\begin{equation}
v(x_*, r_*, \vph_*) = \int \limits_0^{2\pi} \int \limits_0^{\infty}
g({\rm r}_* , {\rm r}) v({\rm r}) r\,dr \,d\vph
\label{continuation}
\end{equation}
(see \cite{Vlasov1995}).

Formulate Green's theorem for the parabolic equation \cite{Sommerfeld1949}. 
Let $\Omega$ be a finite connected domain with 
sectionally smooth boundary $\ptl \Omega$ and an outward normal vector ${\rm n}$. Consider  solutions $v(x,y,z)$, $w(x,y,z)$ of inhomogeneous parabolic equations
\begin{equation}
\label{ineq}
\left(\frac{\ptl}{\ptl x}+ \frac{1}{2ik} \Delta_\perp \right)v(x,r,\varphi) = q(x,r,\varphi),
\end{equation}
\begin{equation}
\label{adjoint}
\left(-\frac{\ptl}{\ptl x}+ \frac{1}{2ik} \Delta_\perp \right)w(x,r,\varphi) = h(x,r,\varphi),
\end{equation}
for some source functions $q$ and $h$. 
Note that the second equation is the adjoint one, i.~e.\  it describes waves traveling in the 
negative $x$-direction.

Introduce vector functions
\begin{equation}
{\bf \rm v}(x,y,z) = \left(\i k v, \frac{\ptl v}{\ptl r},\frac{1}{r}\frac{\ptl v}{\ptl \varphi}\right), \quad {\bf \rm w}(x,y,z) = \left(-\i k w, \frac{\ptl w}{\ptl r},\frac{1}{r}\frac{\ptl w}{\ptl \varphi}\right).
\end{equation}
Here the first element of the vector is the $x$ component, the second is the $r$ component, 
and the third is the $\varphi$ component.
The following relation can be derived using  Gauss--Ostrogradsky theorem:
\begin{equation}
\label{Green}
\int_{\ptl \Omega}[({\bf \rm v}\cdot{\bf \rm n})w-({\bf \rm w}\cdot{\bf \rm n})v]dS = 2 i k \int_\Omega[qw-hv]dV.
\end{equation}
This is the Green's theorem for the parabolic equation.

Note that there are two ways to obtain a solution of an adjoint equation (\ref{adjoint})
from  the initial equation (\ref{ineq}). One can either apply  a complex conjugation or reverse 
the $x$ coordinate.

%%%%%%%%%%%%%%%%%%%%%%%%%%%%%%%%%%%%%%%%%%%%%%%%%%%%%%%%%%%%%%%%%%%%%%%%%%%%

\subsection{Boundary integral equation for the field}
%The equation (\ref{pareq}) describes only waves moving in forward direction due to the first order %$x$-derivative. The following  adjoint parabolic equation describes  backward propagation:
%\begin{equation}
%\label{pareqadj}
%\left(-\frac{\ptl}{\ptl x}+ \frac{1}{2ik}\left(\frac{1}{r}\frac{\ptl}{\ptl r}\left(r\frac{\ptl }%%
%{\ptl r}\right)+ \frac{1}{r^2}\frac{\ptl^2}{\ptl \varphi^2}\right)\right)u = 0.
%\end{equation}
%We will need this equation to state Greens theorem for the parabolic equation.
%Let us state  Green's theorem for the parabolic equation.

\begin{figure}[ht]
\center{\includegraphics[width=6cm]{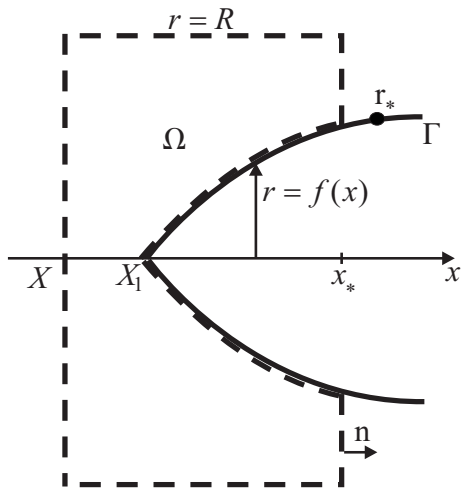}}
\caption{Domain for Green's theorem}
\label{fig01}
\end{figure}

Apply Green's theorem to the domain whose cross-section is shown
in the Fig.~\ref{fig01}. The domain is bounded by the planes
$x = X$, where $X < X_1$, and $x = x_*$ for some $x_* > X_1$, by a cylinder $r = R$ ($R$ is large;
later on the limit $R \to \infty$ is taken),
and by the scatterer surface $\Gamma$ (i.~e.\ $r = f(x)$). Domain 
$\Omega$ has the axial symmetry.    

Let $v$ in (\ref{ineq}) be equal to $u^{\rm sc}$. Obviously, $q\equiv 0$. Let 
$w$ be equal to $g({\rm r}_*, {\rm r})$. Here ${\rm r}_*$ is the source of the field, and it is fixed, while
${\rm r} = (x, r, \vph)$. Note that the source and the observation points are interchanged comparatively to the definition (\ref{Greens}). Thus, the function obeys the adjoint equation 
(\ref{adjoint}) with $h = \delta({\rm r} - {\rm r}_*)$.   

Take the source point for $w$ as
\[
{\rm r}_* = (x_* + \eps , f(x_* + \eps) , \vph_*) 
\]
for arbitrary $\vph_*$ and small $\eps$. The source belongs to $\Gamma$, but it does not belong to 
$\Omega$. Thus, for the Green's theorem $h \equiv 0$ inside $\Omega$.

The right-hand side of (\ref{Green}) is equal to zero.  
The integral in the left can be split into four integrals 
$I_1$, $I_2$, $I_3$, $I_4$, 
related, respectively to the surfaces formed by 
$x = X$, $x = x_*$, $r = R$, and $\Gamma$. 
Integral $I_1$ is equal to zero due to the initial condition (\ref{initcond}). 
Integral $I_3$ turns to zero as $R \to \infty$ due to the radiation condition.
The conical point makes no problem for the integration due to the Meixner's condition 
and Cauchy–-Bunyakovsky inequality.  
Thus, (\ref{Green}) reads as
\[
I_2 + I_4 = 0.
\]
Consider integral $I_2$. According to the structure of (\ref{Green}), it is equal to 
\begin{equation}
\label{I2}
I_2 = 2 i k \int \limits_{f(x_*)}^{\infty} \int \limits_0^{2 \pi}
v(x_* , r, \vph) w(x_* , r, \vph) r d\vph dr.
\end{equation}
Since $w$ is the Green's function with the source located near the plane $x = x_*$, the 
integral can be simplified. 
A direct computation leads to 
\[
\lim_{R \to \infty} I_2 = i k u^{\rm sc}(x_* , f(x_*), \vph_*) + O(\eps) ,
\] 
and
\begin{equation}
\label{I21}
\lim_{\eps \to 0}
\lim_{R \to \infty} I_2 = i k u^{\rm sc}(x_* , f(x_*), \vph_*).
\end{equation}

Now consider integral $I_4$: 
\[
I_4 = \int \limits_{\ptl \Omega \cap \Gamma}
\frac{1}{\sqrt{1 + (\dot f)^2}} \left[ 
v \frac{\ptl w}{\ptl r}
-
w \frac{\ptl v}{\ptl r} 
+
2 i k \dot f v w 
\right]
dS =
\]
\[
\int \limits_{X_1}^{x_*}
\int \limits_0^{2 \pi}
f(x)
 \left(
u^{\rm sc} \bar N [w] - w N[u^{\rm sc}]  
\right)|_{r = f(x)}
d \vph
 dx,
\]
where we introduced the operator
\begin{equation}
\label{barN}
\bar N \equiv \frac{\ptl }{\ptl r} + i k \dot f .
\end{equation}
According to (\ref{boundsc}), 
\begin{equation}
\label{I41}
I_4 = 
\int \limits_{X_1}^{x_*}
\int \limits_0^{2 \pi}
f(x)
 \left(
u^{\rm sc} \bar N [w] + w N[u^{\rm in}]  
\right) \Big|_{r = f(x)}
d \vph
 dx .
\end{equation}
Finally, (\ref{Green}) can be rewritten as 
\[
U^{\rm sc} (x_* , \vph_*) = 
\frac{i}{k} 
\int \limits_0^{2 \pi}
\int \limits_{X_1}^{x_*}
U^{\rm sc}(x,\vph) \bar N  [w] (x,\vph)
f(x)
dx 
d \vph
+ \qquad \qquad \qquad
\]
\begin{equation}
\label{integral_equation_1}
\qquad \qquad \qquad
\frac{i}{k} 
\int \limits_0^{2 \pi}
\int \limits_{X_1}^{x_*}
W(x,\vph) N  [U^{\rm in}] (x, \vph)
f(x)
dx 
d \vph ,
\end{equation} 
where the capital letters correspond to the values on the surface of the scatterer:
\[
U^{\rm sc} (x, \vph) \equiv u^{\rm sc} (x,f(x), \vph),
\]
\[
W (x, \vph) \equiv w (x,f(x), \vph),
\]
\[
\bar N[W] (x, \vph) \equiv \bar N [w] (x,f(x), \vph),
\]
\[
N[W^{\rm in}] (x, \vph) \equiv  N [u^{\rm in}] (x,f(x), \vph).
\]

Equation (\ref{integral_equation_1}) is an inhomogeneous boundary integral equation. Note that 
function $w (x, \vph)$ depends on the variables $x_* , \vph_*$ as on parameters. Thus, 
$\bar N [w]$ plays the role of the kernel of this equation.  

Equation (\ref{integral_equation_1}) can be simplified. For this, apply (\ref{Green})
to 
\[
v({\rm r}) = u^{\rm in} ({\rm r}),
\qquad 
w ({\rm r}) = g({\rm r}_* , {\rm r}) 
\]
in  $\Omega$ shown in Fig.~\ref{fig01}. Also use relation (\ref{continuation}). As the result, get the relation
\[
U^{\rm in} (x_* , \vph_*) = \frac{i}{k} \int \limits_0^{2 \pi}
\int \limits_{X_1}^{x_*}
U^{\rm in}(x,\vph) \bar N  [W] (x,\vph)
f(x)
dx 
d \vph
-
\]
\begin{equation}
\label{integral_equation_2}
\frac{i}{k} 
\int \limits_0^{2 \pi}
\int \limits_{X_1}^{x_*}
W(x,\vph) N  [U^{\rm in}] (x, \vph)
f(x)
\, dx \, d \vph
+ 2 U^{\rm in} (x_*, \vph_*),
\end{equation} 
where 
\[
U^{\rm in} (x, \vph) \equiv u^{\rm in} (x,f(x), \vph).
\]
Finally, introduce the kernel function 
\begin{equation}
\label{kernel}
K(x_* , \vph_* , x, \vph) \equiv
\frac{i f(x)}{k} 
\bar N[W](x, \vph) =
\frac{i f(x)}{k} 
\bar N[g({\rm r}_*, \cdot)](x, \vph) 
\end{equation}
and rewrite (\ref{integral_equation_1}) and (\ref{integral_equation_2})
as a boundary integral equation for the total field on the surface
$U = U^{\rm in} + U^{\rm sc}$:
\begin{equation}
\label{integral_equation_3}
U (x_* ,  \vph_*) = 
\int \limits_0^{2 \pi}
\int \limits_{X_1}^{x_*}
K(x_* , \vph_* , x , \vph)
U(x,\vph) \,
dx \,
d \vph
+ 2 U^{\rm in} (x_*, \vph_*).
\end{equation}

The explicit form of the kernel is as follows: 
\[
K(x_*,\vph_*, x, \vph ) = 
\frac{i k f(x)}{2 \pi}
\left[
\frac{\dot f(x)}{x_*-x}
+
\frac{f(x) - f(x_*) \cos (\vph - \vph_*)}{(x_* - x)^2}
\right]
\times
\qquad \qquad \qquad \qquad 
\]
\begin{equation}
\label{explicit_kernel}
\qquad \qquad \qquad \qquad
\exp \left\{
\frac{i k}{2}
\frac{f^2(x_*) + f^2(x) - 2 f(x_*) f(x) \cos(\vph - \vph_*)}{x_* - x}
\right\}
\end{equation}

Note that equation (\ref{integral_equation_3}) is an equation of Volterra type with respect 
to variable $x$, and is of difference kernel type with respect to the variable~$\vph$.
Both features are of great importance for our consideration.  

Integral equations derived here are close to those obtained in \cite{Hong1967} in ray coordinates.

%%%%%%%%%%%%%%%%%%%%%%%%%%%%%%%%%%%%%%%%%%%%%%%%%%%%%%%%%%%%%
\subsection{Boundary integral equation for the angular modes}

Use the rotational symmetry of the geometry of the problem. 
Let the incident field and the total field (both taken on the surface of 
the scatterer) be represented as the Fourier series: 
\begin{equation}
U^{\rm in} (x, \vph) = 
\sum_{n = - \infty}^{\infty}
U^{\rm in}_n (x) e^{i n \vph},
\qquad 
U(x, \vph) = 
\sum_{n = - \infty}^{\infty}
U_n (x) e^{i n \vph}.
\label{Fourier}
\end{equation}
One can see that functions $U_n(x)$ obey the following integral equations: 
\begin{equation}
U_n(x_*) = \int \limits_{X_1}^{x_*}
K_n (x_* , x) U_n(x) \, dx + 2 U^{\rm in}_n (x_*),
\label{1DIE}
\end{equation}
where $K_n$ are the Fourier components of the difference kernel
\begin{equation}
K_n (x) = 
 \int \limits_0^{2\pi} K (x_* , \vph , x, 0) e^{- i n \vph} d\vph.
\label{1Dkernel}
\end{equation}

Using (\ref{explicit_kernel}), obtain an explicit expression for the kernels
$K_n (x_* , x)$:
\[
K_n(x_*,x) = \frac{ikr}{(x_*-x)^2}\exp \left\{
\frac{i k}{2}
\frac{r_*^2 + r^2}{x_* - x}
\right\}\times
\]
\begin{equation}
\label{Kn}
\left[(-i)^n\left(r+(x_*-x)\dot f(x)\right)J_n\left(\frac{i k r_*r}{x_*-x}\right)  + (-i)^{n+1}r_*\dot J_{n+1}\left(\frac{i k r_*r}{x_*-x}\right)\right],
\end{equation}
 where $r = f(x)$, $r_* = f(x_*)$, 
 $J_n$ is the Bessel functions of the first kind.

%%%%%%%%%%%%%%%%%%%%%%%%%%%%%%%%%%%%%%%%%%%%%%%%%%%%%%%%%%%%%%%%%%%%%%%%%%%%%%%%%
\subsection{Reconstruction of the field outside the surface of the scatterer}

Since there exists Green's theorem for the parabolic equation, one can build formulae
similar to those for Hemholtz equation. Let equations (\ref{1DIE}) be solved, and let the 
total field be known on the surface of the obstacle.  
This field is denoted by $U(x,\vph)$. 
Reconstruct the field everywhere outside the obstacle. 

Apply (\ref{Green}) to the fields 
\[
v({\rm r}) = u({\rm r}), 
\qquad 
w({\rm r}) = g({\rm r}_* , {\rm r}),
\qquad
{\rm r}_* = (x_* , r_*, \vph_*)
\] 
for $r_* > f(x)$. 
The domain for (\ref{Green}) is $X_1 - \eps < x < x_* + \eps$.
Apply relation  (\ref{continuation}). 
As the result, get 
\begin{equation}
u ({\rm r}_*) = \frac{i}{2 k}
\int \limits_0^{2 \pi}
\int \limits_{X_1}^{x_*} 
\bar N[g({\rm r}_* , {\rm r})]
U(x, \vph) f(x) \, dx \, d\vph
+ 
u^{\rm in} ({\rm r}_*),
\label{reconstruction}
\end{equation} 
where operator $\bar N$ acts on the second argument of $g$, and 
${\rm r} = (x, f(x), \vph)$.

Let the scatterer be local, i.~e.\ occupy some domain 
$X_1<x<X_2$, $r < f(x)$, $X_2 < x_*$.
According to (\ref{reconstruction}),  
the scattered field is as follows: 
\begin{equation}
u^{\rm sc} ({\rm r}_*) = \frac{i}{2 k}
\int \limits_0^{2 \pi}
\int \limits_{X_1}^{X_2} 
\bar N[g({\rm r}_* , {\rm r})]
U(x, \vph) f(x) \, dx \, d\vph.
\label{reconstruction_1}
\end{equation} 
Let be  
\[
{\rm r}_* = (L ,  L \theta_*, \vph_*)
\]
for fixed $\theta_*$, $\vph_*$, and large $L$. 
Expand the phase as power series
\[
g({\rm r}_* , {\rm r}) \approx 
\frac{k}{2 \pi i L}
\exp \left\{ \frac{i k L \theta_*^2}{2} \right\}
\exp \{ i k (- \theta_* r \cos (\vph_* - \vph) + x \theta_*^2 / 2) \} 
\]
and transform (\ref{reconstruction_1}) as 
\begin{equation}
u^{{\rm sc}}(L, \theta_* L , \vph_*) \approx
\frac{k}{2\pi i L}
\exp \left\{ \frac{i k L \theta^2}{2} \right\}
T(\theta_* , \vph_*).
\label{reconstruction_2}
\end{equation}
Here $T$ is the directivity (diffraction coefficient): 
\begin{equation}
T(\theta_* , \vph_*) =
\frac{1}{2}
\int \limits_0^{2 \pi}
\int \limits_{X_1}^{X_2} 
(\theta_* \cos(\vph - \vph_*) - \dot f(x)) \times
\qquad \qquad \qquad \qquad
\label{reconstruction_3}
\end{equation}
\[
\exp \{ i k (- \theta_* f(x) \cos (\vph_* - \vph) + x \theta^2 / 2) \} 
U(x, \vph)
f(x)
\,
dx
\,
d\vph.
\]
This integral can be found using the Fourier series representation of $U$ (\ref{Fourier}):
\begin{equation}
T(\theta_*,\varphi_*) =\frac{(-i)^n}{2}\sum\limits_{n=-\infty}^{\infty}
\int \limits _{X_1}^{X_2}\left[\theta_*\dot J_{n+1}(k\theta_* f(x)) -\dot f(x)  J_n(k\theta_* f(x))\right]\times
\end{equation}
\[
\exp\left\{ i k x\theta_*^2/2\right\}U_n(x)f(x)dx.
\]

The optical theorem can be formulated for this case \cite{Korolkov2016e1,Korolkov2016e2}: 
 
\begin{equation}
\int\limits_{0}^{\infty}\int\limits_{0}^{2\pi}|u^{\rm sc}(x , r , \vph)|^2 r dr d\varphi = -2{\rm Re}[T(\theta,0)],
\end{equation}
for some $x>X_2$.
%%%%%%%%%%%%%%%%%%%%%%%%%%%%%%%%%%%%%%%%%%%%%%%%%%%%%%%%%%%%%%%%%%%%%%%%%%%%%%%%%
\section{Diffraction by a cone for axial incidence}

%%%%%%%%%%%%%%%%%%%%%%%%%%%%%%%%%%%%%%%%%%%%%
\subsection{Boundary integral equations for the conical problem}

\begin{figure}[ht]
\center{\includegraphics{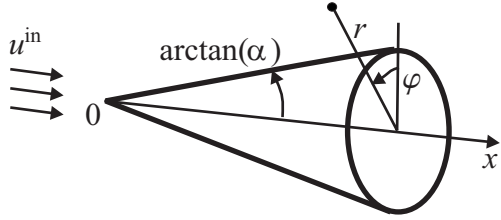}}
\caption{Geometry of the problem}
\label{fig02}
\end{figure}

Consider a cone located at $x > 0$, i.~e.\ $X_1 =0$.  
The profile function is a straight line (\ref{coneshape})
%\begin{equation}
%\label{coneprofile}
%f(x) = \alpha x,
%\end{equation}
where $\alpha$ is the tangent of the angle 
between the axis of the cone (the positive $x$-axis) and the generatrix of the cone. 
Assume that $\alpha\ll 1$.  
The vertex corresponds to the origin (see Fig. \ref{fig02}).

%%%%%%%%%%%%%%%%%
%%%%%%%%%%%%%%%%%

The incident wave travels in the positive $x$-direction. The 
incident wave has the simplest form (\ref{uin}). 
The total field is axially-symmetrical, i.~e.\ 
only the zeroth component is present in the series (\ref{Fourier}):
\[
U (x, \vph) = U_0 (x) \equiv U(x). 
\]

Using (\ref{1DIE}) and (\ref{Kn}) obtain the integral equation
 
\begin{equation}
U(x_*) = \int \limits_0^{x_*} K_0 (x_* , x) U(x) dx + 2.
\label{IEcone}
\end{equation}
The kernel for the axially-symmetrical component of the field in the conical
case can be found from (\ref{Kn}): 
\[
K_0(x_*,x) = \frac{ik \alpha^2 x_*x}{( x_*-x)^2}\exp\left\{\frac{ik \alpha^2}{2}\frac{ x_*^2 + x^2}{ x_*-x}\right\} \times
\qquad \qquad \qquad \qquad
\]
\begin{equation}
\qquad \qquad \qquad \qquad
\left[J_0\left(k \alpha^2\frac{ x_* x}{ x_*-x} \right)+iJ_1\left(k\alpha^2 \frac{ x_* x}{
x_*-x}\right)\right].
\end{equation}

In what follows, we concentrate on solving (\ref{IEcone}). 
 
%%%%%%%%%%%%%%%%%%%%%%%%%%%%%%%%%%%%%%%%%%%%%
\subsection{Solution of the integral equation by Fourier method}
\label{subsec}
Represent $K_0(x_*,x)$ as follows:
\begin{equation}
K_0(x_*,x) = \tau^2 \zeta(\tau)\zeta(\tau_*)^{-1}G(\tau-\tau_*), \quad \tau = x^{-1},\quad \tau_* =x_*^{-1},
\end{equation} 
\begin{equation}
G(\xi) = ik \alpha^2 \xi^{-2}
\exp\left\{ik\alpha^2 \xi ^{-1}\right\}
\left[J_0\left(k \xi ^{-1}\alpha^2\right)+
iJ_1\left(k \xi^{-1}\alpha^2\right)\right],
\end{equation}
\begin{equation}
\label{zeta}
\zeta(\tau) = \frac{k}{\tau}\exp\left\{-ik\tau^{-1}\alpha^2 /2\right\},
\end{equation}
Introduce a new unknown function: 
\begin{equation}
\label{newfun}
V(\tau) = \zeta(\tau)U(1/\tau).
\end{equation}
An integral equation for $V$ has form 
\begin{equation}
\label{intur3}
V(\tau_*) = \int\limits^{\infty}_{\tau_*} G(\tau-\tau_*)V(\tau)d\tau + 
2 \zeta(\tau_*).  
\end{equation}
One can see that this equation has a convolution form, and thus it can be 
solved by Fourier method.

Namely, introduce the Fourier transform 
for any function $p(\tau)$ 
and its inverse as follows:
\begin{equation}
\label{ftrans}
\tilde p(\lambda) =\int\limits^{\infty}_{-\infty}
p(\tau) e^{-i\lambda \tau} d\tau,
\qquad 
 p(\tau) =\frac{1}{2\pi}\int\limits^{\infty}_{-\infty}\tilde p(\lambda)
 e^{i\lambda \tau }d\lambda.
\end{equation}

The equation (\ref{intur3}) takes the following form in the Fourier domain:
\begin{equation}
\label{Fourierint}
\tilde V(\lambda) =  \tilde G(\lambda)\tilde V(\lambda)+ \tilde \zeta(\lambda),
\end{equation}
where
\[
\tilde
G(\lambda) = \int \limits^{\infty}_0 G(\tau) e^{i\lambda \tau} d\tau.
\]

After some algebra (see Appendix \ref{appendix}) obtain the following solution: 

\begin{equation}
\label{cone_sol}
U(x) = \frac{i}{2 k x}
\exp\left\{\frac{ik \alpha^2 x}{2} \right\}
\int \limits^{\infty}_{0}
\frac{\dot J_0(\alpha \sqrt{\lambda})H^{(1)}_0(\alpha \sqrt{\lambda})}{
\dot H^{(1)}_0(\alpha \sqrt{\lambda})}
\exp\left\{\frac{i\lambda}{2kx}\right\} d\lambda + 1.
\end{equation}
The integral in the right is, obviously, the field $U^{\rm sc}$ on the surface of the cone, 
and the unity is the incident field.

Using the reconstruction formula (\ref{reconstruction_1}), one can obtain 
an expression for the scattered field outside the cone: 
\begin{equation}
\label{ufull}
u^{\rm sc}(x,r) =  \frac{i}{2kx}\exp\left\{\frac{ik \alpha^2 x}{2} \right\}
\int \limits^{\infty}_{0}\frac{\dot J_0(\alpha \sqrt{\lambda})H^{(1)}_0(\sqrt{\lambda} r/x)}{
\dot H^{(1)}_0(\alpha \sqrt{\lambda})}\exp\left\{\frac{i \lambda}{2kx}\right\} d\lambda.
\end{equation}
 
Formula (\ref{ufull}) coincides with results \cite{Shanin2012c,Andronov2012ae}. 
This is a bit surprising because these results (especially \cite{Shanin2012c}) 
have been obtained using rather different approaches.

Indeed, the possibility to represent the kernel of the integral equation 
in the difference form has some roots in the mathematical nature of the 
problem. Namely, the parabolic problem we formulated admits separation 
of variables. This consideration  is close  to that of 
\cite{Andronov2012ae}. 

The separation of variables can be done in the coordinates
$(x, \rho , \vph)$ with:
\begin{equation}
\rho =  r / x  .
\end{equation}
Also, introduce new field variable $\hat u$:
\begin{equation}
u(x, \rho x , \vph) = \Xi^{-1}(x, \rho) \, \hat u(x, \rho, \vph),
\qquad
\Xi(x, \rho) \equiv kx \exp\left\{- ikx \rho^2 / 2 \right\}.
\end{equation}
(note that $\zeta(1/x) = \Xi(x,\alpha)$).

Equation (\ref{pareq}) takes the form:
\begin{equation}
\label{pareq3}
\left(\frac{\ptl}{\ptl x}+ \frac{1}{2ikx^2}\left(\frac{1}{\rho}\frac{\ptl}{\ptl \rho}\left(\rho\frac{\ptl }{\ptl \rho}\right) +\frac{1}{\rho^2}\frac{\ptl^2}{\ptl \vph^2}\right) \right)
\hat u = 0.
\end{equation}
Boundary condition (\ref{boundcond}) is
\begin{equation}
\label{boundcond2}
\frac{\ptl \hat u}{\ptl \rho}\at[\Big]{\rho = \alpha} = 0.
\end{equation}
The field $\hat u$ can be represented as a sum of scattered and incident field:
$ \hat u = \hat u^{\rm sc} + \hat u^{\rm in} $.
A general form of solution $\hat u^{\rm sc}$ is as follows:
\begin{equation}
\label{sol1}
\hat u^{\rm sc} = \int C(\lambda) H^{(1)}_0( \sqrt{\lambda} \rho)
\exp\left\{i \frac{\lambda}{2kx}\right\} d\lambda,
\end{equation}
for some unknown coefficient $C(\lambda)$ and contour of integration. 
Here $\lambda$ is the separation constant.

Function $C(\lambda)$ should be determined from the boundary condition (\ref{boundcond2}). 
The procedure is rather straightforward and it yields (\ref{ufull}).

%%%%%%%%%%%%%%%%%%%%%%%%%%%%%%%%%%%%%%%%%%%%%%%%%%%%%%%%%%%%%%%%%%%%%%%%%%%%%%%%%%%%%
\subsection{Elementary analysis of solution (\ref{cone_sol})}

\label{elementary}

The first feature following from (\ref{cone_sol}) (indeed, it can be derived from the 
integral equation itself) is the self-similarity of the solution. The solution 
$u(x)$ depends on a single variable $x$ and on two parameters $\alpha$, $k$, but it can be represented as follows: 
\begin{equation}
\label{cone_sol_1}
U^{\rm sc}(x) = \frac{i}{y}
e^{i y/ 2} 
\int \limits^{\infty}_{0}
\frac{\dot J_0(\varkappa)H^{(1)}_0(\varkappa)}{
\dot H^{(1)}_0(\varkappa)}
\exp\left\{\frac{i\varkappa^2}{2y}\right\} \varkappa \, d\varkappa  ,
\qquad
y = k x \alpha^2.
\end{equation} 
For large $y$ this integral can be estimated by using the standard asymptotic methods. Namely,
the estimation consists of two terms: the stationary phase term corresponding to $\kappa = y$ 
(the estimation shows that it is equal to 1)
and the starting point term corresponding to small $\varkappa$. The  result is as follows:
\begin{equation}
U^{\rm sc} (x) = 1 + P \frac{e^{i y / 2}}{y} + o(y^{-1}), 
\label{cone_sol_2}
\end{equation}
where
\[
P \equiv i \int \limits^{\infty}_{0}
\frac{\dot J_0(\varkappa)H^{(1)}_0(\varkappa)}{
\dot H^{(1)}_0(\varkappa)} \varkappa \, d\varkappa .
\]
These terms admit a clear physical interpretation. The stationary phase term corresponds to the 
reflected wave. Since Neumann boundary conditions are studied, the reflection coefficient 
is equal to~1. Another term corresponds to the diffracted rays traveling along the surface of the cone.  
The amplitude decays as $x^{-1}$. The phase corresponds to a plane wave (in the parabolic approximation) traveling at the ``angle'' $\alpha$ to the $x$-axis.

A more subtle analysis shows that  
\begin{equation}
U^{\rm sc} (0)  = 0 .
\label{cone_sol_3}
\end{equation} 
This relation, indeed, is imposed by the initial condition, but it is not easy to be seen from the integral form of the solution.

%%%%%%%%%%%%%%%%%%%%%%%%%%%%%%%%%%%%%%%%%%%%%%%%%%%%%%%%%%%%%%%%%%%%%%%%%%%%%%%%%%%%%
\subsection{Formula (\ref{ufull}) and the structure of the penumbral zone}

There is a considerable amount of papers solving conical diffraction problems. 
Here we reproduce some concepts from these papers. 

\begin{figure}[ht]
\center{\includegraphics{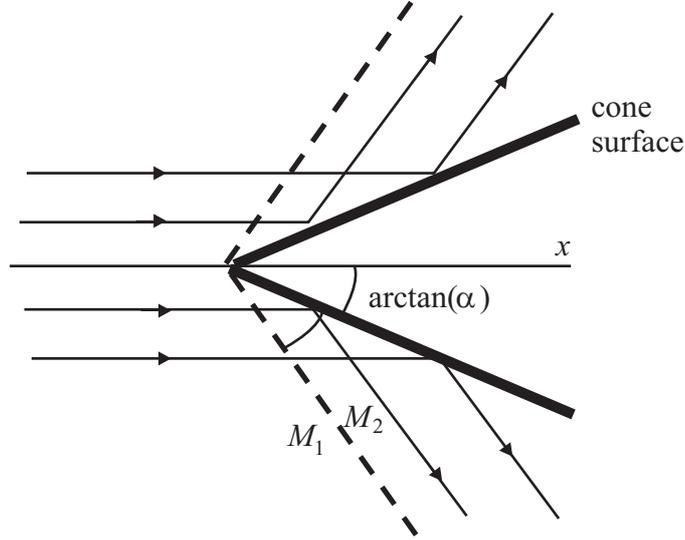}}
\caption{Domains $M_1$ and $M_2$ for a cone}
\label{oasis}
\end{figure}

On the elementary level of understanding the ``ray structure'' of the field 
is shown in Fig.~\ref{oasis}. 
In the case of the axial incidence the whole space outside the cone is illuminated by the 
incident wave. There is a domain (a cone of vertex angle equal to $2 \arctan(\alpha)$)
illuminated by the reflected rays.  This domain is called $M_2$ in \cite{Babich2000}. The 
domain not reached by the reflected rays is called $M_1$ there. In both domains there exists a 
spherical wave diffracted by the tip of the cone. Usually, finding the directivity of this
spherical wave is the main target of the research. Domain $M_1$ is usually called the 
``oasis'' since the Smyshlyaev's formula for this domain \cite{Babich2000} admits contour deformation 
and can be used for efficient computation of the directivity.  

%%%%%%%%%%%%%%%
%%%%%%%%%%%%%%% 

In a narrow domain surrounding the boundary between domains $M_1$ and $M_2$ there is a 
penumbra, since the boundary separates the zone illuminated by the reflected rays and 
the zone not illuminated by them. 
It is well-known (see for example \cite{Borovikov1966en,Babich2000,Babich2004en,Popov2009}) 
that 
the field in this penumbra can be expressed in terms of the 
parabolic cylinder function, namely 
\begin{equation}
D_{-3/2}(z) = \frac{2}{\sqrt{\pi}}\exp\left\{-\frac{1}{2}z^2\right\}
\int \limits_0^{\infty}\exp\left\{-zs-\frac{1}{2}s^2\right\}s^{1/2}ds.
\end{equation}
 
Here our aim is to obtain the similar expression from formula (\ref{ufull}). 
Rewrite (\ref{ufull}) as follows: 
\[
u^{\rm sc}(x, r)  = \frac{i}{y}
e^{i y/ 2} 
\int \limits^{\infty}_{0}
\frac{\dot J_0(\varkappa) }{
\dot H^{(1)}_0(\varkappa)}
H^{(1)}_0 \left( \frac{r }{\alpha x} \varkappa  \right)
\exp\left\{\frac{i\varkappa^2}{2y}\right\} \varkappa \, d\varkappa. 
\]

Again, an estimation of this integral for large $y$ can be obtained 
by using the stationary phase method. 
There are, generally, two terms, one corresponding to the stationary phase point
\[
\varkappa = \left(2 - \frac{r}{\alpha x}  \right) y
\] 
and another corresponding to the end point $\varkappa = 0$. If these points are close to each other, one should 
simplify the integral by taking the leading terms in the Taylor series of the phase and by 
obtaining a ``standard integral'' peculiar to the problem. One can see that the 
case 
\[
2 - \frac{r}{\alpha x} \approx 0
\]  
corresponds to the penumbral zone, and the standard integral in this case is the parabolic 
cylinder function. The field in the penumbra can be written as
\begin{equation}
\label{PC}
u^{\rm sc} \approx -\frac{ie^{i\pi/8}}{(kx)^{1/4}}\sqrt{\frac{x}{r}}\exp\left\{i \frac{k x}{2}\tan^2\alpha -kx\gamma^2 \right\}D_{-3/2}\left(\sqrt{kx}\gamma e^{i3\pi/4}\right),
\end{equation}
where $\gamma = \left(2\alpha- r / x\right)$. The result (\ref{PC}) is in agreement with results of works \cite{Borovikov1966en,Babich2000,Babich2004en,Popov2009}.

%%%%%%%%%%%%%%%%%%%%%%%%%%%%%%%%%%%%%%%%%%%%%
\subsection{Numerical solution of the integral equation by iterations}

Equation (\ref{IEcone}) belongs to the Volterra type. Thus, it can be solved by iterations. 
Namely,
\begin{equation}
U(x) = \sum_{n = 0}^{\infty} U^{(n)} (x),
\label{iterations_1}
\end{equation} 
\begin{equation}
U^{(0)}(x) = 2 U^{\rm in} (x) = 2,
\label{iterations_2}
\end{equation} 
\begin{equation}
U^{(n+1)}(x_*) = \int \limits_0^{x_*} K_0(x_* , x) U^{(n)} (x) \, dx ,
\qquad n > 0.
\label{iterations_3}
\end{equation} 

\begin{figure}[ht]
\center{\includegraphics[width=0.8\linewidth]{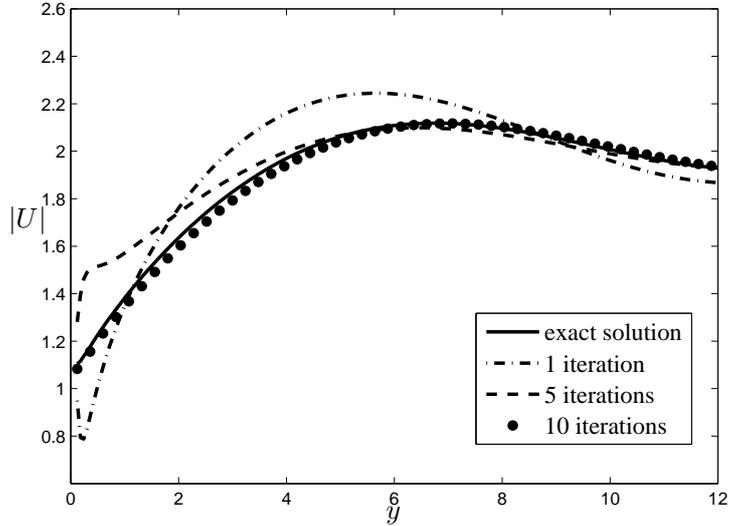}}
\caption{Numerical calculation of the iteration series. Here $y = \alpha^2kx$}
\label{numerics}
\end{figure}

The terms of series can be calculated numerically. The results are shown in Fig.~\ref{numerics}.
One can see that several
 terms (about 10) are enough to model the field 
with a good accuracy.   

%%%%%%%%%%%%%%%%%%%%%%%%%%%%%%%%%%%%%%%%%%%%%
\subsection{Some speculations on the solution by iterations}

\begin{figure}[ht]
\center{\includegraphics{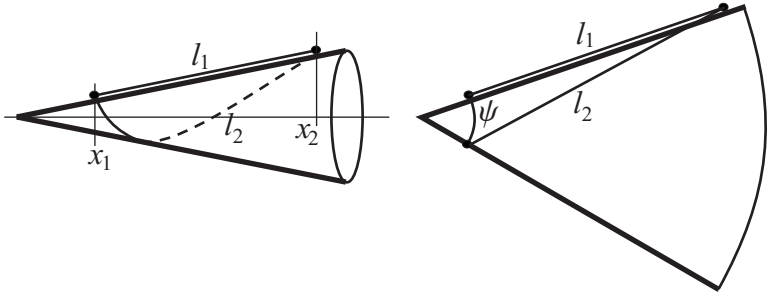}}
\caption{To interpretation of $x_{\rm c}$}
\label{speculations_1}
\end{figure}

In Subsection~\ref{elementary} we noticed that the solution depends only on 
the dimensionless coordinate $y = \alpha^2 k x$. For $y \gg 1$ 
the representation (\ref{cone_sol_2}) is valid. Thus, there exists 
an important spatial scale  $x_{\rm c} = (k \alpha^2)^{-1}$. 
This size can be interpreted as follows. Take two points on the cone surface, 
with coordinates $(x_1, \vph=0)$ and $(x_2 , \vph=0)$. 
Let $x_2$ be very large. Consider two geodesic paths along the cone surface 
shown in Fig.~\ref{speculations_1}, left. The first one, $l_1$, is the shortest way between 
the points, and second one, $l_2$, makes one turn about the axis of the cone. 

%%%%%%%%%%%%%%%%%
%%%%%%%%%%%%%%%%%  

Estimate the difference $\Delta l = l_2 - l_1$. 
For this, draw the paths on the sweep of the cone shown in Fig.~\ref{speculations_1}, right. 
One can see that for large $x_2$
\[
\Delta l \approx (1 - \cos \psi) \sqrt{1 + \alpha^2} x_1,
\]
where $\psi$ is the angle of the sweep:
\[
\psi = \frac{2 \pi \alpha}{\sqrt{1 + \alpha^2}}
\]
One can see that $\Delta l \sim x_1 \alpha^2$. Thus, the domain 
$x_1 \ll x_{\rm c}$ corresponds to the case of $\Delta l$ much smaller than the wavelength, 
and $x_1 \ll x_{\rm c}$ corresponds to $\Delta l$ much bigger than the wavelength. 

One can say that for $x \gg x_{\rm c}$ diffraction process is angularly local (with respect to a remote observation point), 
and for $ x < x_{\rm c}$ the process is not angularly local. For a local process one can use 
a ray approach (similar to that from \cite{Popov2009}). For a non-local process one should better use 
the modal approach from \cite{Andronov2012ae}.   

The iterative approach to the integral equation has a possible interpretation, as we think. 
As it is known, diffraction by a set of edges can be represented in terms of the diffraction 
series or Schwarzschild's series \cite{Schwarzschild1902}.  Each term corresponds to an act of diffraction. This approach usually provides an efficient solution although there are no physical ``acts of diffraction'' (we should note that in \cite{Borovikov1966en} some physical interpretation for the acts of 
diffraction is given). However, this approach has not been extended to 
the Fock's processes of diffraction by curved surfaces.   

We can propose to consider the terms of the series (\ref{iterations_1}) as the results of  successive acts of diffraction. In the previous work \cite{Shanin2016en} we studied diffraction 
by a parabola in a 2D space and found that the the number of terms necessary to describe 
the diffraction process is of order of the size of the shadow zone measured by the Fock's 
scale $ \Delta x \sim (f'')^{-2/3}k^{-1/3}$. I.~e.\ when a creeping wave passes a Fock's zone size, it can be considered as a single diffraction act.

In the current paper we studied diffraction by a cone whose longitudinal curvature is equal to zero. 
We found that several terms (about 10) are enough to describe the wave field on the surface of the cone. This means that the diffraction process happens locally within several $x_{\rm c}$ from the 
tip of the cone, and then there occurs only a ray propagation  along the surface. 

Some estimations of the kernel of the integral equation can be used to understand where 
diffraction happens and how many terms are necessary to describe the process.
\section{Conclusions}
\label{sec:conclusions}
The problem of diffraction by a slender body of rotation in a 3D space is approximated by 
a parabolic formulation in Cartesian coordinates. Using the Green's formula, a boundary integral equation of Hong's type \cite{Hong1967} is derived for the total field on the surface. 
For a thin cone and the axial incidence the equation is solved analytically (by Fourier 
transform) and numerically (by iterations). It is shown that several iterations 
are enough to describe the field with high accuracy.

The work is supported by RSF grant 14-22-00042.

\appendix
\section{Derivation of the solution of  equation (\ref{intur3})} 
\label{appendix}
Equation (\ref{intur3}) takes form (\ref{Fourierint}) in the Fourier domain. Let us find $\tilde \zeta(\lambda)$. It follows from (\ref{zeta},\ref{ftrans}) that
\begin{equation}
\tilde \zeta(\lambda) = \int\limits_{-\infty}^{\infty}\frac{2k}{\tau}\exp\left\{-ik\alpha^2/(2\tau)-i\lambda\tau\right\}d\tau.
\end{equation}

Using a well-known relation 
\begin{equation}
\exp\left\{\frac{ia}{\tau}\right\}=\frac{i\tau}{2}\int\limits_0^\infty J_0(\sqrt{a}t)\exp\left\{-\frac{i\lambda^2\tau}{4}\right\}t dt,
\end{equation}
one can obtain:
\begin{equation}
\tilde \zeta(\lambda) = -4\pi i k J_0\left(\sqrt{2k\alpha^2\lambda}\right),\quad \lambda>0.
\end{equation}
\begin{equation}
\tilde \zeta(\lambda) = 0, \quad \lambda<0.
\end{equation}
Function $\tilde G(\lambda)$ can be calculated in a similar way:
\begin{equation}
\label{G}
\tilde G(\lambda) = 1+ J_0(\sqrt{2\lambda k\alpha^2})\dot H^{(1)}_0(\sqrt{2\lambda k\alpha^2})i\pi\sqrt{2\lambda k\alpha^2}.
\end{equation}

Thus the solution of equation (\ref{intur3}) is given by the following formula: 
\begin{equation}
\label{prelimsol}
\tilde V(\tau) = \frac{1}{\pi \alpha}\int\limits_{0}^{\infty}\frac{\sqrt{\lambda}}{\dot H^{(1)}_0\left(\alpha\sqrt{\lambda}\right)}\exp\left\{\frac{i\lambda}{2kx}\right\}d\lambda.
\end{equation}
Representation (\ref{cone_sol}) follows from (\ref{prelimsol}) (see \cite{Andronov2012ae}).

%\bibliographystyle{siamplain}
%\bibliography{references_thin_cone}

%%%%%%%%%%%%%%%%%%
%%%%%%%%%%%%%%%%%%%%%%%%%%%%

\ifx\undefined\BibEmph\def\BibEmph#1{#1}\else\fi
\ifx\undefined\href\def\href#1#2{#2}\else\fi
\ifx\undefined\url\def\url#1{\texttt{#1}}\else\fi
\ifx\undefined\urlprefix\def\urlprefix{URL: }\else\fi
\ifx\undefined\BibUrl\def\BibUrl#1{\urlprefix\url{#1}}\else\fi
\ifx\undefined\BibUrlDate\long\def\BibUrlDate#1{({%
\cyr\cyrd\cyra\cyrt\cyra\ %
\cyro\cyrb\cyrr\cyra\cyrshch\cyre\cyrn\cyri\cyrya}: #1)}\else\fi
\ifx\undefined\BibAnnote\long\def\BibAnnote#1{#1}\else\fi

\end{document}